\documentclass[10pt]{easychair}
\usepackage{graphicx}
\usepackage{titlesec}
\usepackage{indentfirst,csquotes}

\usepackage{amssymb,amsthm,amsmath}
\usepackage{xcolor,paralist,hyperref,titlesec,fancyhdr,etoolbox}

\titleformat{\section}[display]{\normalfont\huge\bfseries\centering}{\centering\chaptertitlename\thechapter}{10pt}{\Large}
\titlespacing*{\section}{0pt}{0ex}{0ex}

\hypersetup{ colorlinks=true, linkcolor=black, filecolor=black, urlcolor=black }

\usepackage{lipsum}

\begin{document}
\title{A Naïve Encoding of Russell's Paradox in Type Theory} 
\author{Qu Zhuoyuan}
\institute{Nagoya University}
\maketitle

\titlerunning{A Naïve Encoding of Russell's Paradox in Type Theory}
\authorrunning{QU Zhuoyuan}

\let\thefootnote\relax

\begin{abstract}
Russell's paradox is the most easily understandable way to illustrate the inconsistency of naïve set theory. This note proposes a direct encoding of Russell's paradox with type-in-type universe, sigma types, and either extensional identity or intensional identity with the uniqueness of identity proofs (UIP).
\end{abstract} 

\bigskip

\titleformat{\section}
  {\normalfont\Large\bfseries\raggedright} 
  {\thesection}                            
  {1em}                                    
  {}                                       

\section{Background}

It is well-known that a "type in type" system can construct Girard's paradox\cite{HurGirard}, breaking consistency. The paradox is the type-theoretic analogue of Burali-Forti paradox. But there are not many attempts on Russell's paradox. Coquand\cite{CoqTree} used W-types (general inductive types) to make a tree that corresponds to Russell's paradoxical set. It is a quite old paper, for which 1Lab project has a nice formalization in Agda\cite{1lab}, using an inductive type instead of the raw W-type. In this note, we will show a way to formalize Russell's paradox in a type-in-type system with sigma-types and identity types with UIP.

\section{Russell's paradox in set theory}

In a naïve set theory, Russell constructed a set where all of its elements are not members of the element itself, written as $R := \{ x \in V | x \not\in x \}$, where $V$ is the set of all sets including $V$ itself ($ V \in V$). Then one may ask whether $R \not\in R$. It is true, since if $R \in R$ we should have $R \not \in R$ by the definition of $R$. However, $R \in R$ is also true, since we already have $R \not \in R$, which means that $R$ follows the property of $R$, so it should be in $R$. And that leads to a contradiction.

A common misunderstanding is that the key reason of Russell's paradox is that all the objects in set theory are sets. For example, in set theory $0 := \{\}$, $1 := \{\{\}\}$, so we have a fact that $0 \in 1$, which is meaningless in mathematics. Such theorems are called junk theorems. In Russell's paradoxical set, the property $x \not \in x$ is a meaningless proposition. In type theory (including Russell's type theory), objects are separated into types; it might seem that one cannot construct Russell's paradoxical set.

\section{Constructing Russell's paradox in type theory}
\def\U{\mathcal{U}}
\def\P{\mathcal{P}}
We assume that we are in a type theory with a single universe $\U$ such that $\U : \U$. We define a type representing subsets of all types.

\[
\begin{split}
&V : \U \\
&V := \sum_{A : \U} \P A
\end{split}
\]

Here $\P A := A \to \U$, as the type of predicates on $A$. We will use the following notation,

\[
    \{x : A \ | \ P(x) \} := (A, P) : V
\]

where $P : \P A$.
As examples, we have such definitions,

\[
\begin{split}
    &\mathbb{N} := \{ n : \text{Nat} \ | \ \top \} \\ 
    &\mathbb{Z}^+ := \{ n : \text{Int} \ |\ n > 0 \}
\end{split}
\]

Before we construct the paradoxical set, we need the "element of" operation. Let us try,
\[
\begin{split}
&(\in_A) : A \to V \to \U \\
&a \in_A \{ x : B \ |\  P(x) \} := P(a)
\end{split}
\]
Here $a : A$ but $P : B \to \U$, so we require an operation to transport $x : A$ to $x : B$.

Thus we define,

\def\coe{\text{coe}}
\def\Id{\text{Id}}
\def\refl{\text{refl}}
\def\coeq{\text{coe-eq}}

\[
\begin{split}
&\coe_{AB} : \Id_\U(A, B) \to A \to B\\
&\coe_{AB}\ (\refl, \ a) := a  \\
&\\
&\coeq_A : \prod_{(x : A)} \prod_{h : \Id(A,A)} \Id_A(a, \coe(h, a)) \\ 
&\coeq_A\ (x, \refl) := \refl\\
&\\
&(\in_A) : A \to V \to \U \\
&a \in_A \{ x : B \ |\  P(x) \} := \sum_{h : \Id(A, B)} P(\coe(h, a))
\end{split}
\]

Note that the definition of \coe\ and \coeq\ used dependent pattern matching, \coe\ can be easily defined with the axiom J, but \coeq\ requires not only the axiom J but also the axiom K, or equivalently, UIP.

If we are using extensional type theory\cite{nlab:extensional_type_theory}, we do not even need to define $\coe$, as $a : A$ can be directly transported to $a : B$ under the assumption $h : \Id(A,B)$, as in \cite{hjboom} \cite{ruspar}.

We write $(\in) := (\in_V)$ for short. This means that sets can appear in both sides of $(\in)$. This is exactly what we need in the paradoxical set.

After that we can construct the paradoxical set, 

\[
R := \{x : V \ | \ x \not \in x \}
\]

where $x \not \in x := x \in x \to \bot$.

Thus we have,

\def\la{\text{lemma1}}
\def\lb{\text{lemma2}}
\def\subst{\text{subst}}

\[
\begin{split}
    &\la : R \not \in R \\
    &\la\ (h, p) := p (\subst_{\lambda x. x\in x}(\coe(\coeq(R, h)), (h,p))
\end{split}
\]

Where for $P : A \to \U$, $\subst_P : \Id_A(x, y) \to P(x) \to P(y)$ is defined by,
\[
\begin{split}
    &\subst_P\ (\refl, p) := p
\end{split}
\]

We can illustrate the proof in English. Assuming that we have $H : R \in R$, to show $\bot$, we first split $H$ (since $R \in R$ is a sigma-type) into $h : \Id(V, V)$ and $p : \coe(h, R) \not \in \coe(h, R)$. Since we have $\coeq(R, h) : \Id(R, \coe(h, R))$, we can use it to rewrite the type of $p$ and get $R \not \in R$, which contradicts $H : R \in R$.

With $\la$, constructing $R \in R$ is easy,
\[
\begin{split}
    &\lb : R \in R \\
    &\lb:= (\refl, \la)
\end{split}
\]

\begin{section}{Difference from other paradox constructions}

In \cite{coquand23}, Coquand shows that if one can have the following structure, then some of the paradoxes can be constructed.

\def\intro{\text{intro}}
\def\match{\text{match}}
\[
\begin{split}
    &V : \U \\ 
    &\intro : \P^2 V \to V \\ 
    &\match  : V \to \P^2 V \\
    \text{A judgemental equality: } &\match \ \circ \ \intro = \P^2 \ (\intro \ \circ \ \match)
\end{split}
\]

$\P^2$ should be considered as a functor here. Hurkens' paradox, as a simplification of Girard's paradox, also fits this structure. The idea behind this construction is similar to Russell's idea, but uses a pair of sets that are not mutually included in any other sets, rather than a single set.

If we have the following, similar but simpler structure, we can directly derive Russell's paradox from it.

\[
\begin{split}
    &V : \U \\ 
    &\intro : \P V \to V \\ 
    &\match  : V \to \P V \\
    &\beta_1 : \prod_{u : \P V} \prod_{v : V} \match (\intro(u), v) \to u(v) \\
    &\beta_2 : \prod_{u : \P V} \prod_{v : V} u(v) \to \match (\intro(u), v) 
\end{split}
\]

We can define a paradox as Russell's instead of Girard's, if it can fit into this structure. The construction described in this note is Russell's by taking 
\[
\begin{split}
    &\intro(u) = \{\ x : V\ |\ u(x)\ \}\\ 
    &\match(u,v) = v \in u
\end{split}
\]

Coquand's paradox of trees \cite{CoqTree} can also fit into this structure. See \cite{mine}.

\end{section}

\begin{section}{Conclusion}

This note illustrates that it is possible to naïvely encode Russell's paradox in type theory. However, the construction strongly depends on UIP, hence it will not be able to formalize in homotopy type theory\cite{hottbook} where UIP cannot be assumed. 

In this construction, we use pi, sigma, and identity types to give an instance of Russell's paradox; compared to Coquand's paradox of trees, we avoid using inductive types, but additionally use axiom K here. According to Naïm Camille Favier\cite{ncf}, we can avoid the usage of sigma-types by converting all sigma types to pi types, which also can fit into the structure in the previous section.

Compared to those paradoxes with only pi types, we avoid the occurrence of $\P^2 \U$, which makes the construction much easier to understand.

The whole construction has been formalized in Agda prover and Rocq prover with universe checking off\cite{mine}.

\

In \cite{CoqTree}, Thierry Coquand cited a communication in the TYPES Forum\cite{ruspar}, which is no longer available. We thank Thierry Coquand for finding us H.J. Boom's construction in ETT \cite{hjboom}.

\end{section}

\
\
\
\bibliographystyle{IEEEtran}
\bibliography{bib}

\end{document}